\documentclass[anz]{cupaus}
\usepackage{verbatim}

\begin{document}

\verso{Series solution of Laplace problems}
\recto{Lloyd N. Trefethen}

\title{SERIES SOLUTION OF LAPLACE PROBLEMS}

\cauthormark 
\author[1]{LLOYD N. TREFETHEN}

\address[1]{Mathematical Institute, University of Oxford,
Oxford OX2 6GG, UK\email[1]{trefethen@maths.ox.ac.uk}}

\pages{1}{6}
      
\begin{abstract}
At the ANZIAM conference in Hobart in February, 2018, there were
several talks on the solution of Laplace problems in multiply
connected domains by means of conformal mapping.  It appears to
be not widely known that such problems can also be solved by
the elementary method of series expansions with
coefficients determined by least-squares fitting on the boundary.
(These are not convergent series; the coefficients depend on
the degree of the approximation.)  Here we give a tutorial
introduction to this method, which converges at an exponential
rate if the boundary data are sufficiently well-behaved.
The mathematical foundations go back to Runge in 1885
and Walsh in 1929.
One of our examples involves an approximate Cantor set with up
to 2048 components.
\end{abstract}

\keywords[\textit{Keywords and phrases}]{Green function, conformal mapping,
harmonic measure, Laplace problem, series expansion, least-squares,
Cantor set}

\maketitle

\section{Introduction}\label{sec-intro}  

I am a card-carrying conformal mapper~\cite{SC}, but it is my
view that conformal mapping is usually not the best strategy
for solving Laplace problems in multiply connected regions.
The trouble is that the conformal mapping problem is typically as
difficult as the original Laplace problem, or more so, because
it may require resolution of geometric issues that are absent
from the original problem.  The aim of this paper is to give
a tutorial introduction to the simple alternative of solving
Laplace problems by series expansion with least-squares matching
of boundary data.  Nothing here is mathematically new, though
these methods are not as well known as they might be.

An important tool in connection with multiply connected
conformal mapping is the Schottky--Klein prime
function~\cite{ckgn,crowdy1,green,nasser}, which gives
considerable insight into the structure of such a map.  If the
prime function is known for a particular domain, then it can be
used to solve certain Laplace problems~\cite{crowdy08,crowdy12}.
However, the prime function itself can only be computed
numerically, and in fact, one of the best methods for computing
it makes use of the same kind of series expansions reviewed
here~\cite{crowdy,ckgn}.  All in all, though the idea of solving
multiply connected Laplace problems by conformal mapping is an
old one---a book on such matters was published thirty years ago
by Prosnak~\cite{prosnak}---this approach is computationally unnecessary.

The aim of this paper is to show how series expansion methods
can be used.  It is structured in a tutorial fashion, treating
a succession of problems, with MATLAB code listings given in
the appendix.  These codes, which are descendants of codes in
the unpublished essay~\cite{tda}, are an essential part of the
presentation, and they include details that are not spelled out
in the text.

For concreteness, most of the problems take the form of the
computation of a Green function with its logarithmic singularity
at the origin.  This amounts to the computation of the harmonic
measure~\cite{gm,ransford}, that is, the function that determines the
probability density that a particle undergoing Brownian motion
from the origin will first encounter the boundary at each
particular point.  We also give a few examples illustrating how
other kinds of Laplace problems, involving non-constant boundary
conditions for example, can be treated by the same approach.

Many people have discovered versions of these ideas over the
years, and I doubt anyone has a comprehensive understanding of the
literature.  Some references are given in Section~\ref{sec-hist},
starting with Runge's theorem of 1885.

\section{Green function for a disk}\label{sec-disk}  

We begin with a problem that could be solved analytically.
What is the Green function $u(z)$ in the exterior of a disk of radius $r$
centered at position $z=c$ in the complex plane?\ \ Specifically,
we seek a function $u$ satisfying
\begin{equation}
\Delta u = 0,
\quad u(z) = 0 \hbox{ for } |z-c\kern .3pt|=r,
\quad u(z) \sim \log|z| \hbox{ as } z\to 0
\label{prob}
\end{equation}
in the region of the complex plane exterior to $z=0$ and the
circle $|z-c\kern .3pt|=r< |c|$, with regular behavior at $z=\infty$,
i.e., $u(z)\to u_\infty^{}$ as $z\to\infty$ for some constant
$u_\infty^{}$.  Throughout this paper, we regard the plane as
either real or complex according to convenience.  Thus the function
$\log|z|$ in (\ref{prob}), for example, is a real function in the
$x$-$y$ plane that has been expressed for simplicity in terms of
the complex variable $z= x+iy$.

To solve (\ref{prob}), we approximate $u$ by a series expansion
\begin{equation}
u(z) = \log|z| - \log|z-c\kern .3pt| + C + \sum_{k=1}^N \left[a_k^{}\kern 1pt
\hbox{Re\kern 1pt} ((z-c)^{-k}) + b_k^{}\kern 1pt
\hbox{Im\kern 1pt }((z-c)^{-k})\right],
\label{expdisk1}
\end{equation}
which could be written in real form as
\begin{equation}
u(z) = \log |z| - \log r + C + \sum_{k=1}^N r^{-k} \kern 1pt [a_k^{}
\cos(k\theta) - b_k^{}\sin(k\theta) ]
\label{expdisk1real}
\end{equation}
with $z-c = r\kern .7pt e^{i\theta}$.  The coefficients $a_k^{}$
and $b_k^{}$ are chosen to satisfy $u(z) = 0$ as nearly as
possible in a least-squares sense at ${\tt \kern 2pt npts}\gg N$
sample points along the boundary.  The term $\log|z-c\kern .3pt
|$ is needed to make $u$ regular as $z\to\infty$.  (One could
alternatively use $\log|z-\tilde c\kern .3pt |$ with $\tilde
c$ equal to a different point in the interior of the disk.)
Note that only negative powers of $z-c$ appear in the series,
because the domain is unbounded.  Since the functions in play are
analytic and the boundaries are smooth, geometric convergence
occurs as a function of $N$, and a modest value like $N=10$
is more than enough for plotting accuracy.

\begin{figure}
\centering
\vspace{1em}
\includegraphics[scale=.65]{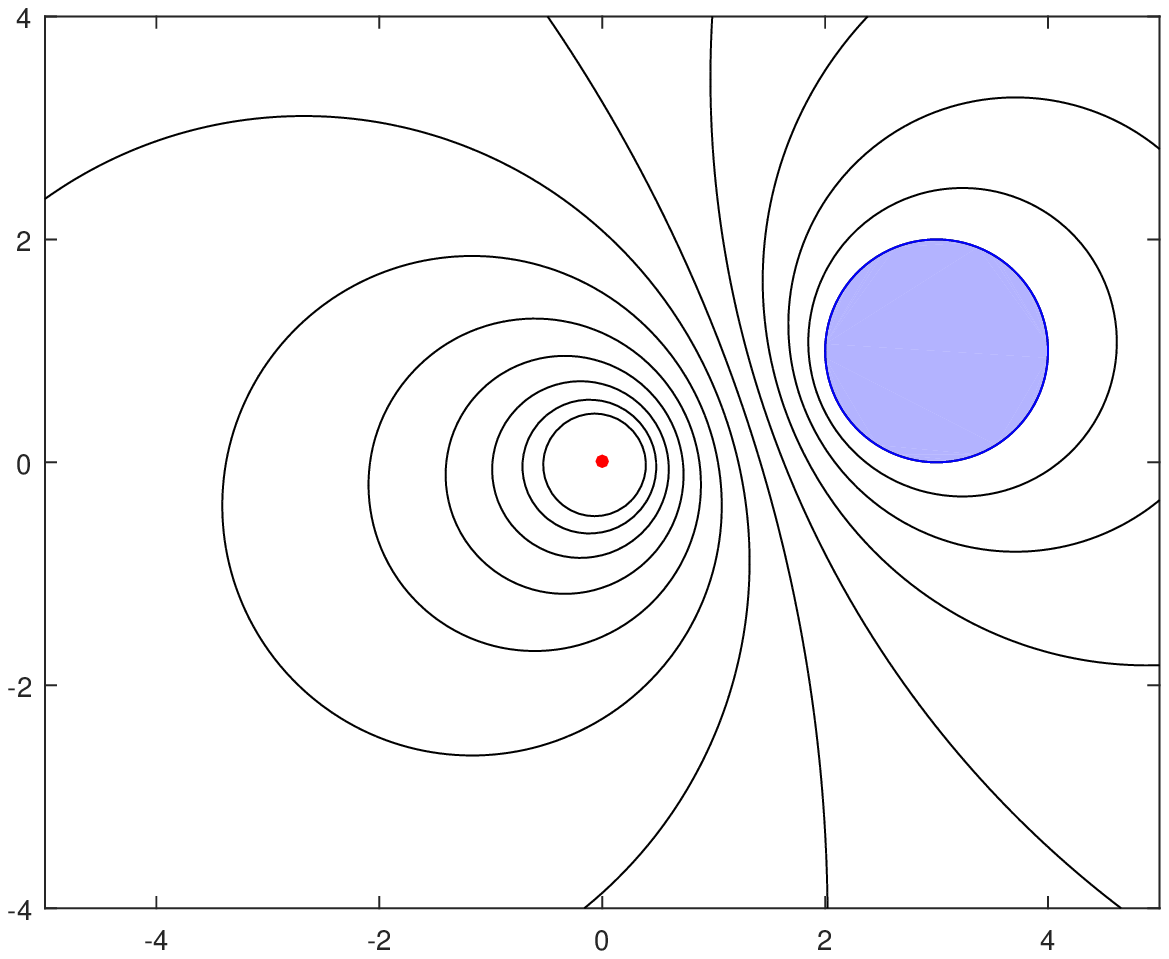}
\caption{Green function $u(z)$ outside a disk computed
with {\tt disk1}.  The level curves
are regularly spaced between $0$ on the boundary of the disk and
$-\infty$ at the singularity at $z=0$. \label{fig-disk1}}
\end{figure}

The MATLAB code {\tt disk1} computes the solution of (\ref{prob})
and plots contour lines for the particular choices $c=3+i$
and $r=1$.  Here as in all our codes, we strive to make the
computational first half compact but legible, following the
philosophy of~``Ten Digit Algorithms''~\cite{tda}.  Here is
that first half of {\tt disk1}:

{\scriptsize
\begin{verbatim}
      % disk1.m    Green function exterior to a disk
      %
      % We seek the function u(z) that is zero for |z-c| = r and harmonic outside this
      % circle, including at z=infty, except with u(z)~log|z| as z->0.  u is expanded as
      %
      %   u(z) = log|z| - log|z-c| + a(1)
      %              + SUM_{k=1}^N a(2k)*real((z-c)^(-k)) + a(2k+1)*imag((z-c)^(-k)) .
       
      c = 3+1i; r = 1;                           % center and radius of disk
      N = 10; npts = 3*N;                        % no. expansion terms and sample pts
      z = c + r*exp(2i*pi*(1:npts)'/npts);       % sample points
      rhs = -log(abs(z)) + log(abs(z-c));        % right-hand side 
      A = ones(npts,2*N+1);
      for k = 1:N                                % set up least-squares matrix
        A(:,2*k  ) = real((z-c).^-k); 
        A(:,2*k+1) = imag((z-c).^-k); 
      end
      a = A\rhs;                                 % solve least-squares problem
\end{verbatim}
\par}

\noindent
The second
half of each code, devoted to plotting, has been made extremely
compact and is not very legible; see the appendix.

The level curves computed by {\tt disk1} are shown in
Figure~\ref{fig-disk1}.  As a check of accuracy, we find by
computing with larger values of $N$ that the arbitrary value
$u(2) \approx -0.5893274981708$ is accurate to 4 digits with
$N=4$, to 7 digits with $N=8$, and to 10 digits with $N=12$.
(Rigorous estimates can be based on the maximum principle; see
Section~\ref{sec-details}.  All three computations, including
plotting, take much less than a second on a laptop.

It is worth emphasizing that the series in (\ref{expdisk1}), like
all our series, is not made from Taylor or Laurent coefficients,
which would be independent of $N$ and could achieve convergence
as $N\to\infty$ only in restricted regions of the plane.
The coefficients in (\ref{expdisk1}) depend on $N$, and
they allow convergence irrespective of the shape of any region of
analyticity or harmonicity.  This situation is analogous to
the well-known effect that on the real interval $[-1,1]$,
an analytic function $f$ can be approximated by degree $N$
polynomials with geometric convergence as $N\to\infty$, but
these cannot be Taylor polynomials unless $f$ is analytic in
a disk of radius ${>}\kern 1pt 1$~\cite{atap}.

\begin{figure}
\centering
\vspace{1em}
\includegraphics[scale=.65]{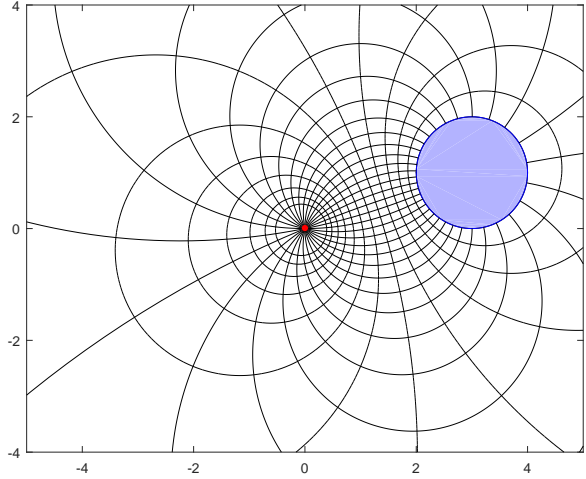}
\caption{Repetition of Figure~\ref{fig-disk1} but now
with an orthogonal system of curves also included (streamlines),
computed with {\tt disk1ode} by solving the ODE (\ref{ODE}).}
\label{fig-disk1ode}
\end{figure}

The curves in Figure~\ref{fig-disk1}
are equipotentials $u(z) = \hbox{const}$, and as always
with harmonic functions, it is also interesting to 
examine an orthogonal system of
curves corresponding to streamlines.  Among other benefits,
these give a visual indication of harmonic measure
along the boundary of the circle.\footnote{The probabilistic
interpretation of a streamline is as follows.  Diffusion
is governed by the heat equation $u_t = \Delta u$, and a harmonic
function $u$ corresponds to a steady state
satisfying the Laplace equation $\Delta u = 0$.
A streamline is a curve across which there is no net diffusion in
this steady state.
If we think of the diffusion process as resulting from particles
moving along Brownian paths,
this means that the rate of particles crossing
from left to right at each point along the streamline
equals the rate crossing from right to left.}
For this simple example, we
could compute the orthogonal system just as level curves of the
imaginary part of the same expansion.  For multiply connected
domains, however, that approach leads to headaches of tracking
branches of the complex logarithm.  A simpler alternative
is to compute the orthogonal curves by solving an ODE to climb
the gradient,
\begin{equation}
\frac{dz}{dt} = \frac{\nabla u}{\| \nabla u\|},
\label{ODE}
\end{equation}
starting from points close to the singularity at $z=0$ and
equally spaced around it.
The gradient can be computed by taking advantage of the complex
variables interpretation to note that if $f$ is a complex analytic function,
then the gradient of its real part is given (in the form of 
a vector represented as a complex number) by
\begin{equation}
\nabla [\hbox{\kern 1pt Re\kern 1pt} f(z)] = \overline{f'(z)}.
\label{gengrad}
\end{equation}
In particular, for real $a,b,d$ we have
\begin{equation}
u(z) = d\log|z| ~~ \Longrightarrow ~~ \nabla u(z) = d/\kern 1.5pt \overline{z}
\end{equation}
and
\begin{equation}
u(z) = a\kern 1.5pt \hbox{Re\kern 1pt }(z^{-k}) +
b\kern 1.5pt \hbox{Im\kern 1pt }(z^{-k})  ~~
\Longrightarrow ~~ \nabla u(z) = -k\kern .7pt (a+bi\kern .7pt )
\kern 1pt  /\kern 1.5pt \overline{z}^{\kern 2pt k+1}.
\end{equation}
Our code for climbing the gradients of Figure~\ref{fig-disk1}
is called {\tt disk1ode}, and the results are shown in
Figure~\ref{fig-disk1ode}.  ODE methods for tracking contours have
been used by Gautschi and Waldvogel, as presented in chapter~25
of~\cite{gander}, by Weideman (for climbing poles of solutions
of Painlev\'e equations, unpublished), and no doubt by others.

Careful readers may note that to solve the ODE (\ref{ODE})
numerically, the codes in the appendix call Matlab {\tt ode23}
rather than, as one might expect, {\tt ode45}.  Both choices
work well for regions bounded by disks, but in the case of slits,
we find that {\tt ode45}, with its larger time steps, more often
introduces anomalies as a contour approaches a slit because of
hopping over the slit and making erroneous choices of branch.

\section{Green function for several disks}\label{sec-disks}  

\begin{figure}
\centering
\vspace{1em}
\includegraphics[scale=.65]{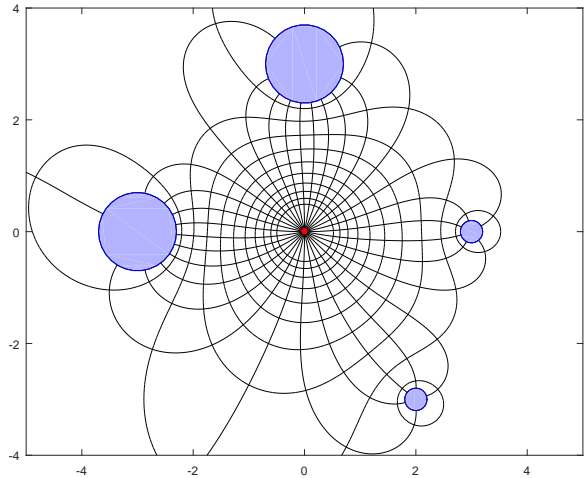}
\caption{An analogous computation with several disks (code {\tt disks}).
Note that even very small disks may have a large influence on the solution.}
\label{fig-disks}
\end{figure}

Our next example is the same as before, except that instead of
one disk defined by $|z-c\kern .3pt|=r$, we have $J$ disks defined by
$|z-c_j^{}| = r_j^{}$, $1\le j \le J$.  The series approximation is 
a generalization of (\ref{expdisk1}),
\begin{equation}
u(z) = \log|z| + C +
\sum_{j=1}^J \left\{ d_j^{} \log|z-c_j^{}| + \sum_{k=1}^N 
\left[ a_{jk}^{}\kern 1pt
\hbox{Re\kern 1pt} ((z-c_j^{})^{-k}) + b_{jk}^{}\kern 1pt
\hbox{Im\kern 1pt }((z-c_j^{})^{-k})\right] \right\},
\label{expdisks}
\end{equation}
together with the condition
\begin{equation}
\sum_{j=1}^J d_j = -1 .
\label{dcond}
\end{equation}
The only new issue that arises in moving from (\ref{expdisk1})
to (\ref{expdisks}) is that the coefficients $d_j^{}$ are
now unknowns, which are determined by imposing the condition
(\ref{dcond}) to ensure nonsingular behavior at $z=\infty$.
We do this by appending one more row to the matrix $A$ and
one more entry to the right-hand side vector,
though one could equally well remove one of the variables
such as $d_1^{}$.  The code is called {\tt disks}.

\begin{figure}
\centering
\vspace{1em}
\includegraphics[scale=.65]{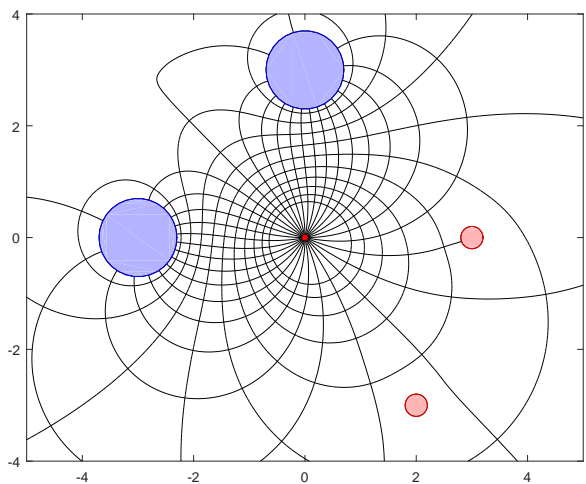}
\caption{Modification in which the two small disks take boundary
values $u=-1$ rather than $u=0$.  This diminishes their influence
greatly.}
\label{fig-disksdiff}
\end{figure}

Figure~\ref{fig-disks} shows the equipotentials and streamlines for a
configuration with four disks, two large and two small.  It is interesting
to record the coefficients $d_1^{},\dots,d_4^{}$, which are the negatives of
the harmonic measures of each disk:
\begin{displaymath}
d_1^{} = -0.342977, ~~
d_2^{} = -0.315902, ~~
d_3^{} = -0.182248, ~~
d_4^{} = -0.158873.
\end{displaymath}
The disks are numbered in clockwise order beginning with the big
one on the left, and the diminishing sizes of these coefficients
correspond to the diminishing proportions of streamlines hitting
each disk.  (We list the numbers to six digits, although with
$N=10$ they have already converged to 13 digits of accuracy.)
Note that although the small disks are much smaller than the
big ones, their harmonic measures are only somewhat smaller.
This is a familiar effect in potential theory, or equivalently
in probability theory: the influence of a structure of radius
$r$ diminishes only at a rate $O(1/\log r)$ as $r\to 0$.
See ~\cite{faraday}, \cite{ransford}, and the Cantor set example
of Section~\ref{sec-cantor}.

Figure~\ref{fig-disksdiff} shows the result of a variant computation.  Here,
instead of $u(z) = 0$ as the boundary condition on all four disks, we require
$u(z) = 0$ on the two larger disks and
$u(z) = -1$ on the two smaller ones (a straightforward change
in the least-squares problem, not listed in the appendix).
The equipotentials and streamlines change
considerably, and the numbers $d_j^{}$ take quite different values:
\begin{displaymath}
d_1^{} = -0.493167, ~~
d_2^{} = -0.492543, ~~
d_3^{} = -0.006514, ~~
d_4^{} = -0.007775.
\end{displaymath}
The tiny values of $d_3^{}$ and $d_4^{}$ show that the small disks
now have little influence on the solution.
This is because the values of $u$ specified on their boundaries are not
much different from the values $u$ would have taken in their
absence.\footnote{Bengt Fornberg points out in
an email: ``If you put $u(z) = +1$ rather
than $-1$ on them, in spite of being small, they would put up a
massive field barrier.  I guess this is the idea behind
the old vacuum tube.''}

\section{Green function for one or several slits}\label{sec-slits}  

\begin{figure}
\centering
\vspace{1em}
\includegraphics[scale=.65]{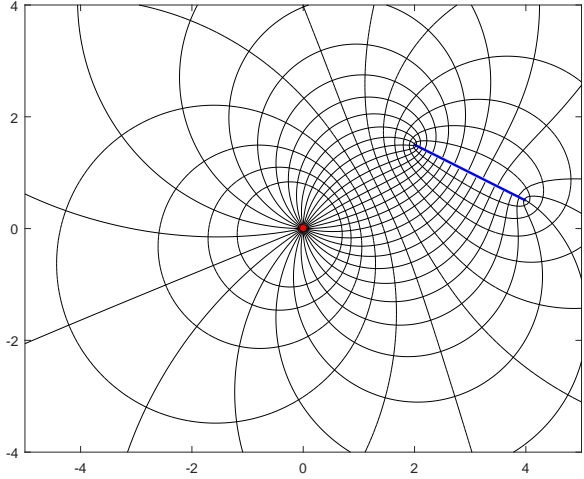}
\caption{Green function $u(z)$ outside a slit computed
with {\tt slit1}. The same method is used as in
Figure~\ref{fig-disk1}, except with basis functions
transplanted from the exterior of a disk to the exterior of the slit.
Computing the difference in imaginary part from one end of the slit to the
other along each side (or equivalently, integrating the
gradient) shows that
58.2625\% of the harmonic measure falls on the side facing
the origin.
\label{fig-slit1}}
\end{figure}

If $c$ and $r\ne 0$ are complex numbers,
then $c + r\kern .7pt [-1,1]$ is a complex
interval that we call a {\em slit,} and the exterior of the unit
disk in the $w$-plane is conformally mapped to the exterior of the slit
by the function 
\begin{equation}
z = c + r\kern .7pt (w + w^{-1})/2.
\label{jouk}
\end{equation}
The inverse map can be written for $z\not\in c + r\kern .7pt [-1,1]$ by
\begin{equation}
w = z_c + s(z_c)\sqrt{z_c^2-1/2}, \quad z_c = (z-c)/r,
\label{jouk2}
\end{equation}
where $s(z_c)$ takes the value $+1$ for $z_c$ in the open right
half-plane or the positive imaginary axis and $-1$ for $z_c$
in the open left half-plane or the negative imaginary axis.
(These choices are designed to work with the standard branch of
the square root.)  For $z\in c+r[-1,1]$, i.e., $z_c\in [-1,1]$,
there are two values of $w$, and we avoid evaluating $w(z)$
for such values.

By transplanting the powers $w^{-k}$, these maps give us good bases for
series expansions in regions with slits.
Equation (\ref{expdisk1}) becomes
\begin{equation}
u(z) = \log|z| - \log|w| + C + \sum_{k=1}^N \left[a_k^{}\kern 1pt
\hbox{Re\kern 1pt} (w^{-k}) + b_k^{}\kern 1pt
\hbox{Im\kern 1pt }(w^{-k})\right],
\label{expslit1}
\end{equation}
and the code {\tt slit1} is an analogue of {\tt disk1}.
Figure \ref{fig-slit1} shows the plot produced by this code,
which uses the choices $c=3$, $d = 1-0.5i$.

\begin{figure}
\centering
\vspace{1em}
\includegraphics[scale=.65]{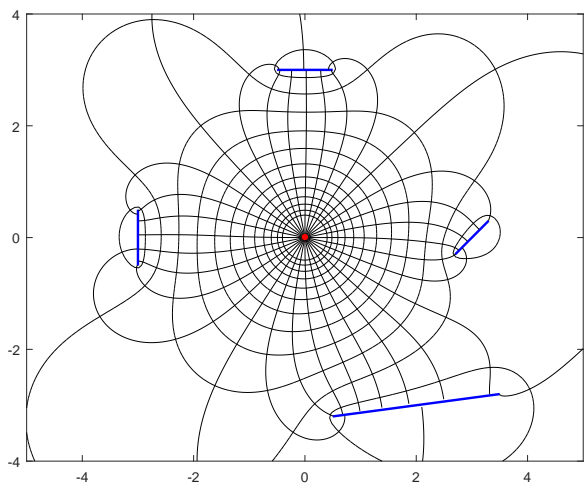}
\caption{Green function $u(z)$ outside several slits, computed
with {\tt slits}.
\label{fig-slits}}
\end{figure}

To treat a region with several slits, we can use this analogue
of (\ref{expdisks}),
\begin{equation}
u(z) = \log|z| + C +
\sum_{j=1}^J \left\{ d_j^{} \log|w_j^{}(z)| + \sum_{k=1}^N 
\left[ a_{jk}^{}\kern 1pt
\hbox{Re\kern 1pt} (w_j^{}(z)^{-k}) + b_{jk}^{}\kern 1pt
\hbox{Im\kern 1pt} (w_j^{}(z)^{-k})\right] \right\},
\label{expslits}
\end{equation}
where $w_j^{}(z)$ denotes the map (\ref{jouk2}) with parameters
$c_j^{}$ and $r_j^{}$.  Without further discussion, we present
in Figure~\ref{fig-slits} the result of executing {\tt slits}.

\section{Approximations to the Cantor set}\label{sec-cantor}
\begin{figure}[p]
\centering
\vspace{1em}
\includegraphics[scale=1.02]{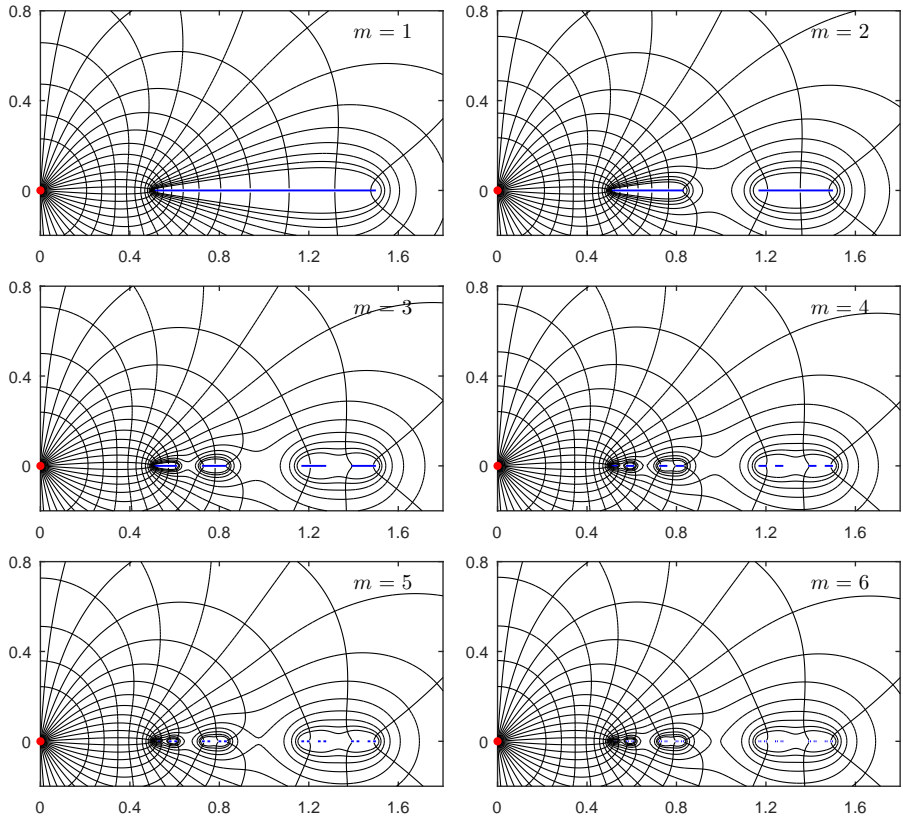}
\caption{Green function $u(z)$ for finite approximations to the 
Cantor set with $2^m$ slits, $m = 1,\dots, 6$, computed
with {\tt cantor}.  The solution is left-right and up-down symmetric, so
only a portion of the right half-plane is shown.  Because of
the $O(1/\log r)$ effect associated with structures of size
$r \ll 1$, the field lines quickly settle down
to nearly their final form.  The computation of these six
images takes a total of about two seconds on a 2016 laptop.
\label{fig-cantor}}
\end{figure}

\begin{figure}[p]
\centering
\vspace{1em}
\includegraphics[scale=1.02]{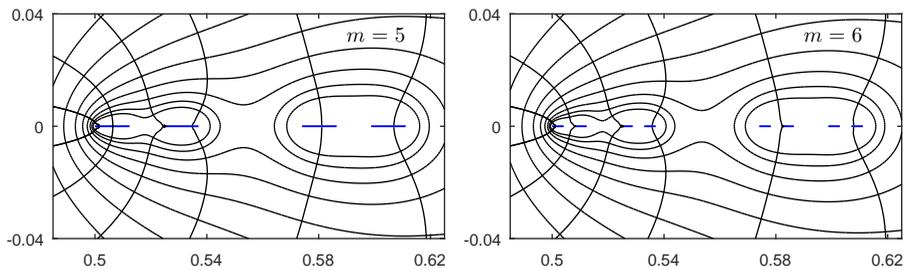}
\caption{Closeups of the last two images of Figure~\ref{fig-cantor}.
\label{fig-closeup}}
\end{figure}

The Cantor set is a canonical example of a fractal, a
subset of the real axis
with infinitely many pieces, zero length, and nonzero capacity.
Figure~\ref{fig-cantor} illustrates that Green functions for
finite approximations to a Cantor set are readily computed
by series expansions.  The plots show that the potential
fields associated with coarse approximations to the Cantor
set quickly settle down to close to the limiting form, a
consequence of the ``familiar effect'' mentioned in the penultimate
paragraph of Section~\ref{sec-disks}.  Figure~\ref{fig-closeup} further
illustrates this with closeups for $m=5$ and $m=6$.

Starting with the interval $[-1.5,1.5]$ of length 3, these
domains are constructed by removing the middle third of each
remaining piece $m$ times, where $m$ is a positive integer.
To compute the fields of Figure~\ref{fig-cantor}, we could
have used the code {\tt slits} as written, but there are
many symmetries to be taken advantage of.  Because of the real
symmetry, the imaginary terms can be removed from the expansions;
also it is only necessary to sample on one side of each slit.
Because of the even symmetry, one can combine terms in pairs
from the left- and right half-planes, and use sample points only
in one of the half-planes.  And, of course, since the solution
computed is left-right and top-down symmetric, one need only
compute field lines in a quadrant.  A further benefit is that
for $m\ge 3$, the slits are so short and well separated that
a small value of $N$ suffices for 6-digit accuracy or more;
we take $N = \max\{2,6-m\}$.  The code {\tt cantor} used for
Figure~\ref{fig-cantor} exploits these effects, and for plots
with $m = 4,5,\dots,10,$ the times required on a 2016 laptop are
about $0.2$, $0.3$, $0.4$, $0.8$, $2.4$, $12$, and $80$ seconds,
respectively.  These geometries have $32, 64, \dots, 2048$ slits.

The negatives of the coefficients $d_j$ of the calculation
give us the harmonic measures of each slit.  From inside out,
here are the harmonic measures of the $2^{m-1}$ slits in the
right half-plane for the level $m$ finite Cantor sets with $m =
1,2,3,4.$  All these figures are believed to be accurate to the
six digits listed: \smallskip

$m=1$:~  1/2
\smallskip

$m=2$:~  0.367776, 0.132224
\smallskip

$m=3$:~  0.253289, 0.111676, 0.066706, 0.068329
\smallskip

$m=4$:~  0.162063, 0.088794, 0.058116, 0.054538,\hfill\break
\indent \kern 2in 0.038156, 0.029363, 0.029460, 0.039509
\smallskip

\noindent Analysis of such numbers confirms that once the slits
in a finite Cantor set are small, not much changes except very
close to the slits.  For example, suppose one adds up the first
halves of the harmonic measures as listed above to determine
the total harmonic measure of the half of the slits closer to
the origin.  For $m=4,5,\dots, 10$ one finds that the measures
add up to about $0.367776$, $0.364965$, $0.363512$, $0.362773$,
$0.362397$, $0.362205$, $0.362107$.  These sums converge to a
limit of about $0.362007$.  Doubling this figure to account for
the left half-plane, we see that the inner slits of a Cantor
set correspond to about 72.4\% of the harmonic measure.

Other numerical computations related to the Cantor set can be
found in~\cite{lsn} and~\cite{ranros}.

\section{Other Laplace problems}\label{sec-other} The examples
we have shown all involve an unbounded domain, a logarithmic
singularity, boundaries consisting of circles or slits,
and piecewise constant boundary conditions. We now indicate
how series methods can be applied for problems without these
features.

{\em Variable coefficients.}  Here it is simply a matter of
sampling the boundary condition in the obvious fashion.  If the
boundary condition and the boundary itself are smooth,
one can expect rapid convergence.

{\em Bounded domains.} The presence of an outer boundary makes it
necessary to include positive-degree terms in a series expansion.
For example, suppose we have a Laplace problem in an annulus
$r_1^{} < |z| < r_2^{}$ with boundary data prescribed on the
inner and outer circles.  Here we would need a series involving
the real parts of both positive and negative (and zero) powers
of $z$, as well as $\log(z)$.  Note that this is analogous
to the Laurent series appropriate for an analytic function in
an annulus, except for the inclusion of the $\log (z)$ term,
which is not relevant to analytic functions since its imaginary
part is not single-valued.

\begin{figure}[t]
\centering
\vspace{2em}
\includegraphics[scale=.7]{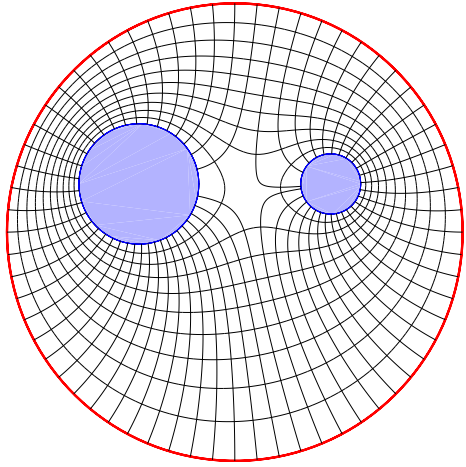}
\caption{Solution of a Laplace problem in a bounded region.
The boundary conditions are $u=0$ on the outer circle and $u=1$ on
the inner ones.
\label{fig-bounded}}
\end{figure}

{\em More general boundary shapes.}  
Smooth boundary components can be treated by the same methods
we have demonstrated.  For a nonsmooth boundary component that
is not simply a slit, one must expect slow convergence unless
special steps are taken.  One approach is to include additional
terms capturing the local behavior of singularities at corners,
as is explored for the Helmholtz rather than Laplace equation
in~\cite{bt,bt2,fhm}.  The important point to note in series
methods is that singularities need only be approximated locally;
the boundary matching will take care of the global connections.
This is in contrast to methods based on conformal mapping, which
require one to resolve global relationships of singularities
before progress is made on the Laplace problem.

Figure~\ref{fig-bounded}, produced by the code {\tt bounded},
illustrates the solution of a Laplace problem in a bounded domain
with no log singularities.  The weights $d_2^{}$ and $d_3^{}$
on the terms $\log|z-z_2^{}|$ and $\log|z-z_3^{}|$ come out as
$-0.840513$ and $-0.450685$.

\section{Mathematical foundations and history}\label{sec-hist}

The ideas underlying series expansion methods date to Runge's
theorem of 1885~\cite{rungethm}.\footnote{Runge is a hero of mine,
who moved from pure to applied mathematics during his career and
late in life was appointed at G\"ottingen effectively as one of
the first professors of numerical analysis in the world.  He was
born 99 years before me to the day.} Runge showed that if $f$ is
an analytic function on a simply connected compact set $K$ in the
complex plane, it can be approximated arbitrarily closely on $K$
by polynomials.  More generally, if $K$ is multiply connected,
say with $J$ holes, then $f$ can be approximated on $K$ by
rational functions with poles at arbitrarily prescribed points
in each of the holes.  Here is how the
theorem is stated by Gaier~\cite{gaier}, with $K^C$ denoting
the complement of $K$ in the complex plane~$\mathbb{C}$.

\medskip

{\narrower\em
\noindent {\bf Theorem 2 {\rm (Runge 1885)}}.  Suppose $K$ is compact in $\mathbb{C}$ and $f$
is analytic on $K$; further, let $\varepsilon>0$.  Then there exists a rational function 
$R$ with poles in $K^C$ such that
\begin{displaymath}
|f(z) - R(z) | < \varepsilon \quad (z\in K) .
\end{displaymath}
\par}

\medskip

\noindent In general these polynomial and rational approximations
cannot come from convergent series, such as Taylor polynomials,
but must have coefficients that vary with the degree $N$.
It is enough for the poles of $R$ to be restricted to lie in
a fixed set of arbitrarily chosen points $c_1^{},\dots , c_J^{}$,
one in each hole.

Our concern here has been harmonic functions 
rather than analytic ones.  The two
settings are almost the same, since the real part of an
analytic function is harmonic, and conversely, a harmonic function
$u$ has a harmonic conjugate $v$ such that $u + iv$ is analytic.
The complication is that in a multiply connected domain, the
function $v$ will in general be multiple-valued, taking different
values along different paths around the holes (or equivalently,
single-valued on a Riemann surface).  This situation is
pinned down in the following statement
from a beautiful paper by Sheldon Axler called ``Harmonic functions
from a complex analysis viewpoint''~\cite{axler}:

\medskip

{\narrower\em
\noindent {\sc Logarithmic Conjugation Theorem}.  Suppose $\Omega$
is a finitely connected region, with $K_1^{},\dots,K_N^{}$
denoting the bounded components
of the complement of $\Omega.$
For each $j$, let $a_j^{}$ be a point in
$K_j^{}$.  If $u$ is a real valued harmonic function on $\Omega,$ then there
exist an analytic function $f$ on $\Omega$ and real numbers $c_1^{},\dots,c_N^{}$
such that
\begin{displaymath}
u(z) = \hbox{Re\kern 1pt} f(z) + c_1^{} \log|z-a_1^{}| + \cdots
+ c_N^{} \log |z-a_N^{}|
\end{displaymath}
for every $z$ in $\Omega.$
\par}

\medskip

\noindent This theorem leads to the generalization of Runge's
theorem from analytic to harmonic functions, which consists
essentially of adding a logarithmic term $d_j
\log|z-c_j|$ corresponding to each hole, exactly as in
our codes.  This was first spelled out
out in a paper by Walsh in 1929~\cite{walsh29} (which Axler
describes as ``the only place I have been able to find the
Logarithmic Conjugation Theorem written down with a proof'').

Thus the foundations of series methods for Laplace problems
are a century old.  It is curious that it is so difficult
to find mention of these foundations in the literature
of these methods.  What happened over the years is that the
mathematicians moved on to ever more refined theorems, always
seeking the weakest possible regularity assumptions (and rarely
mentioning harmonic functions).  Building on work by
Keldysh in the late 1930s, a landmark of such developments
was Mergelyan's theorem of 1951, which allows $f$ to be just
continuous on $K$ and analytic in its complex interior, if any,
rather than analytic on all of $K$~\cite{gaier}.  Such results
are important, but they are quite technical, and remote
from most applications.

Walsh and his student Curtiss had an interest in the idea of
constructing harmonic approximations based on boundary fitting, though they
did not do computations~\cite{curtiss}.  Unfortunately,
like many authors,
they focussed on interpolation rather than least-squares.
Interpolation formulations, with
their square matrices, bring serious challenges of distribution
of points and convergence.  (Curtiss writes: ``The all-important
issue in the practical applications will certainly be the correct
choice of the interpolation points on the boundary $C$ of the
region.'')  When we switch to least-squares, the matrices become
rectangular and such difficulties go away so long as one samples
in sufficiently many points.  Possibly the first advocate of
using least-squares for boundary-matching with Laplace problems
was Cleve Moler in a Stanford technical report of 1969~\cite{moler69}.

Meanwhile, with little influence from the literature of theoretical
mathematics, series expansions were being used in innumerable
applications.  An early example was Lord Rayleigh, in 1892, who
used series to solve Laplace problems in a plane with a periodic
array of circular holes removed~\cite{rayleigh}; a related
method for an electromagnetic problem was applied by Z\'avi\v ska
in 1913~\cite{zaviska}.  The Soviet Union in the mid-20th century
was a center of great expertise in the solution of boundary value
problems for linear partial differential equations (PDEs), with
leading names including Gakhov, Kantorovich, Krylov, Mikhlin,
and Mushkelishvili.  In 1936 Kantorovich and Krylov published
their monograph {\em Methods for the Approximate Solution
of Partial Differential Equations,} which evolved into {\em
Approximate Methods of Higher Analysis,} the latest edition in
1960~\cite{kk}.  Though the style is not very computational by
modern standards, the book has a great deal of material about
series expansions for PDEs.

In the computer era,
series methods have been applied not only for the
Laplace equation but also for the Helmholtz~\cite{bb,bt,bt2},
biharmonic~\cite{bourot2,price}, and Maxwell equations,
among others, and this is the basis for methods with
names including Trefftz's method, the method of particular
solutions~\cite{bergman,bt,bt2,fhm,henrici}, point-matching,
the method of plane waves, and the method of fundamental
solutions~\cite{bb,ogata}, also known as the charge
simulation method~\cite{amano,malik,ssw}, Concerning the
Laplace equation, papers making use of series methods include
\cite{faraday}, \cite{crowdy}, \cite{ddep} (section~4),
and \cite{finn}.  These stem from my own methods reported
in~\cite{tda}, but there are undoubtedly other unrelated
publications.  One such that I am aware of is the very
nice paper by Rostand~\cite{rostand}.

Series methods have much the same flavour as the numerical
solution of integral equations~\cite{nmia}, a very large
subject---though they are not the same, as
explained in~\cite{akt}.  Integral equations methods tend to
take over when the problems feature complicated geometries
or non-smooth data.  For sufficiently complex problems, one
finds oneself in the worlds of $\cal H$-matrices and the Fast
Multipole Method, which enable extremely
complicated complicated geometries to be treated with
nearly linear complexity~\cite{greengard,hack,helsing,nasser1}.
The starting point of the Fast Multipole Method is a recursive
use of just the kinds of series we have described.  Although the
present paper deals with two-dimensional problems, all these
methods apply in 3D too~\cite{cheng}.

\section{Numerical considerations}\label{sec-details}  

The aim of this paper has been to present series methods, not analyze
them.  The starting point for analysis will be complex or harmonic
approximation theory of the kind mentioned in the last section.
Without going into any details, we mention some issues that must
arise in a deeper treatment of series methods.

{\em Boundary resolution.}  The great convenience of a
least-squares formulation is that it permits one to put plenty
of points along the boundary, thereby bypassing subtle questions
related to interpolation.  But how big must the ratio ${\tt
npts}/N$ be to achieve satisfactory accuracy?   And can this
come with guarantees?  A mathematical analysis will require
attention to matters of polynomial or trigonometric interpolation
related to the Runge phenomenon~\cite{runge01} and the study of
Lebesgue constants, and sharp answers will depend on regions
of analyticity of the solutions being computed~\cite{atap,bb}.

{\em Convergence rate.}   As $N\to\infty$, the convergence will
normally be exponential, at a rate $O(C^{-N})$ for some $C>1$,
if the boundaries and the boundary data are analytic.  The value
of $C$ will depend on matters of analytic continuation, as
above~\cite{bb}.  Only in simple cases would one expect
to analyze such rates analytically.

{\em Scaling.}  We have used powers $(z-c)^{-k}$ for all disks,
regardless of their radius.  In problems with widely varying
scales, however, it may be advantageous to scale such terms
relative to the radius.  This is related to column
scaling of the matrix~$A$.

{\em Complexity.}  Series methods rely on solving
least-squares matrix problems, and the cubic complexity
will eventually show up if problems with enough components
are considered.  In such cases one may turn to hierarchical
approaches such as the Fast Multipole Method.  As a rule of thumb,
it seems clear that one should do this (in 2D) for problems with
thousands of components, but perhaps not if there are merely
hundreds of components.

{\em Nearly-touching boundary components.}  In some problems two
holes, say, may nearly touch, and convergence of a series method
will degrade.  In simple cases one may solve the problem by using
a point other than the center for a local expansion (or relatedly,
a local M\"obius transformation).  In more complicated problems,
this is one of the situations where hierarchical methods prove
their power~\cite{chengg1,chengg2,helsing}.

{\em Three dimensions.} Nothing that we have done is restricted
to two dimensions, apart from the use of a complex variable
for convenience~\cite{cheng}.  In 3D, nevertheless, there is no
denying that computations tend to require more human and computer
effort, and the role of simple series methods is perhaps smaller
than in 2D.

{\em Maximum principle.}  One of the attractions about boundary
matching for elliptic problems is that the maximum principle
allows easy {\em a posteriori} analysis of accuracy.  So long
as the boundary data are closely matched, as can be verified by
sampling finely, one can be assured that the computed solution
is accurate.  See for example~\cite{fhm,moler69}.

\section{Conclusion}\label{sec-concl}
The elementary solution methods we have presented are not
always familiar to those who might find them useful.  They tend
to be eclipsed by more general, more powerful tools, which are
necessary for sufficiently complicated problems but rely on much
more machinery.  An explanation for this situation may lie in
the disparate rates of development of mathematics, algorithms,
computers, and software.  Fifty years ago, thanks to Runge and
Walsh and Keldysh and others, the mathematical
basis of series methods was already in place and known (if not in
mathematical detail) to the early computer-era numerical analysts.
And so one finds Hockney, for example, publishing a paper in
1964 about series expansions for a Laplace problem involving ``a
round hole in a square peg''~\cite{hockney}.  But it took time
for computers to grow powerful enough, and software convenient
enough, to enable such methods to be so easily used as we have
demonstrated.  By the time that had happened, researchers had
come to focus on more advanced tools.

As our Cantor set example shows, elementary methods can work
well even for regions with hundreds of holes.

\acks Many people have given me helpful advice along the way,
including Kaname Amano, Alex Barnett, Darren Crowdy, Tom De Lillo,
Toby Driscoll, Bengt Fornberg, Andrew Gibbs, Abinand Gopal,
Christopher Green, Leslie Greengard, Dave Hewett, J\"org Liesen,
Cleve Moler, Yuji Nakatsukasa, Mohamed Nasser, Tom Ransford,
Elias Wegert, and Andr\'e Weideman.  This article was written
during an extremely agreeable sabbatical visit to the Laboratoire
de l'Informatique du Parall\'elisme at ENS Lyon hosted
by Nicolas Brisebarre, Jean-Michel Muller, and Bruno Salvy.

\medskip

\section*{Appendix.  MATLAB codes}

\smallskip

{\scriptsize

\begin{verbatim}
function disk1   % Green function exterior to a disk, equipotential curves only
%
% We seek the function u(z) that is zero for |z-c| = r and harmonic outside this
% circle, including at z=infty, except with u(z)~log|z| as z->0.  u is expanded as
%
%   u(z) = log|z| - log|z-c| + a(1)
%              + SUM_{k=1}^N a(2k)*real((z-c)^(-k)) + a(2k+1)*imag((z-c)^(-k)) .
 
c = 3+1i; r = 1;                           % center and radius of disk
N = 10; npts = 3*N;                        % no. expansion terms and sample pts
z = c + r*exp(2i*pi*(1:npts)'/npts);       % sample points
rhs = -log(abs(z)) + log(abs(z-c));        % right-hand side 
A = ones(npts,2*N+1);
for k = 1:N                                % set up least-squares matrix
   A(:,2*k  ) = real((z-c).^-k); 
   A(:,2*k+1) = imag((z-c).^-k); 
end
a = A\rhs;                                 % solve least-squares problem

% Contour plot
x = linspace(-5,5,145); y = linspace(-4,4,115); [xx,yy] = meshgrid(x,y);
zz = xx+1i*yy; uu = disk1fun(zz); z = c + r*exp(pi*1i*(-50:50)'/50);
fill(real(z),imag(z),[.7 .7 1]), hold on
plot(z,'b'), plot(0,0,'.r','markersize',13)
levels = -3:.25:-.25; contour(xx,yy,uu,levels,'k'), axis equal, axis([-5 5 -4 4])
set(gca,'xtick',-4:2:4,'ytick',-4:2:4,'fontsize',8), hold off, print -depsc disk1

function u = disk1fun(z)
u = log(abs(z)) - log(abs(z-c)) + a(1);
for k = 1:N, u = u + a(2*k)*real((z-c).^(-k)) +a(2*k+1)*imag((z-c).^(-k)); end
u(abs(z-c)<=r) = NaN;
end

end
\end{verbatim}
\smallskip

\begin{verbatim}
function disk1ode   % Green function exterior to a disk
%
% We seek the function u(z) that is zero for |z-c| = r and harmonic outside
% this circle, including at z=infty, except with u(z)~log|z| as z->0.
% u is expanded as
%
%   u(z) = log|z| - log|z-c| + a(1)
%              + SUM_{k=1}^N a(2k)*real((z-c)^(-k)) + a(2k+1)*imag((z-c)^(-k)) .
    
c = 3+1i; r = 1;                           % center and radius of disk
N = 10; npts = 3*N;                        % no. expansion terms and sample pts
z = c + r*exp(2i*pi*(1:npts)'/npts);       % sample points
rhs = -log(abs(z)) + log(abs(z-c));        % right-hand side 
A = ones(npts,2*N+1);
for k = 1:N                                % set up least-squares matrix
   A(:,2*k  ) = real((z-c).^-k); 
   A(:,2*k+1) = imag((z-c).^-k); 
end
a = A\rhs;                                 % solve least-squares problem
    
% Contour plot
x = linspace(-5,5,145); y = linspace(-4,4,115); [xx,yy] = meshgrid(x,y); zz = xx+1i*yy;
uu = disk1fun(zz); z = c + r*exp(pi*1i*(-50:50)'/50); fill(real(z),imag(z),[.7 .7 1])
hold on, plot(z,'b'), plot(0,0,'.r','markersize',6), levels = -3:.25:-.25;
contour(xx,yy,uu,levels,'k'), axis equal, axis([-5 5 -4 4])
set(gca,'xtick',-4:2:4,'ytick',-4:2:4,'fontsize',8), op = odeset('events',@event);
for t = pi*(1:32)/16, z0 = .01*exp(1i*t); sol = ode45(@dzdt,[0 500],z0,op);
plot(deval(sol,linspace(0,max(sol.x),300)),'k'), end
plot(0,0,'.r','markersize',13), hold off, print -depsc disk1ode

function u = disk1fun(z)
u = log(abs(z)) - log(abs(z-c)) + a(1);
for k = 1:N, u = u + a(2*k)*real((z-c).^(-k)) + a(2*k+1)*imag((z-c).^(-k)); end
u(abs(z-c)<=r) = NaN;
end

function g = dzdt(t,z)
g = 1./conj(z) - 1./conj(z-c);
for k = 1:N, g = g - k*(a(2*k)+1i*a(2*k+1))./conj(z-c).^(k+1); end
g = g./abs(g);
end

function [val,isterm,dir] = event(t,z)
dir = 0; isterm = 1; val = abs(z-c)-r;
end

end
\end{verbatim}
\smallskip

\begin{verbatim}
function disks   % Green function exterior to several disks
%
% We seek the function u(z) that is zero for |z-c(j)| = r(j), j = 1,...,J,
% and harmonic outside these circles, including at z=infty, except with 
% u(z)~log|z| as z->0.  u is expanded as
%
%   u(z) = log|z| + C + SUM_{j=1}^J {d(j)*log|z-c(j)|
%            + SUM_{k=1}^N [a(j,k)*real((z-z(j))^-k + b(j,k)*imag((z-z(j))^-k]}
%
% with SUM d(j) = -1; all these coefficients are collected in the vector X.
% The unknowns determined by linear least-squares are C, d(1),...,d(J),
% {a(j,k)}, {b(j,k)}.  A has dimensions 1+J*npts by 1+J*(2N+1).
     
c = [-3 3i 3 2-3i];                         % centers
r = [.7 .7 .2 .2]; J = 4;                   % radii
N = 10; npts = 3*N;                         % no. expansion terms and sample pts
circ = exp(2i*pi*(1:npts)'/npts);           % roots of unity
z = []; for j = 1:J
   z = [z; c(j)+r(j)*circ]; end             % sample points on the bndry
A = ones(size(z));                          % constant term
for j = 1:J
   A = [A log(abs(z-c(j)))];                % logarithmic terms
   for k = 1:N                              % set up least-squares matrix
      zck = (z-c(j)).^(-k);
      A = [A real(zck) imag(zck)];          % algebraic terms
   end
end
A = [A; zeros(1,1+J*(2*N+1))];
A(end,2:(2*N+1):end) = 1;
rhs = [-log(abs(z)); -1];
X = A\rhs;                                  % solve least-squares problem
C = X(1); X(1) = [];                        % extract results
d = X(1:2*N+1:end), X(1:2*N+1:end) = [];
a = X(1:2:end); b = X(2:2:end);
    
% Contour plot
x = linspace(-5,5,145); y = linspace(-4,4,115); [xx,yy] = meshgrid(x,y); zz = xx+1i*yy;
uu = disksfun(zz); z = exp(2i*pi*(0:30)'/30);
for j = 1:J, disk = c(j)+r(j)*z; fill(real(disk),imag(disk),[.7 .7 1])
hold on, plot(disk,'b'), end, levels = -1.9:.2:-.1;
contour(xx,yy,uu,levels,'k'), axis equal, axis([-5 5 -4 4])
set(gca,'xtick',-4:2:4,'ytick',-4:2:4,'fontsize',8), op = odeset('events',@disksevent);
for t = pi*(1:32)/16, z0 = .01*exp(1i*t); sol = ode45(@dzdt,[0 100],z0,op);
plot(deval(sol,linspace(0,max(sol.x),300)),'k'), end
plot(0,0,'.r','markersize',13), hold off, print -depsc disks

function u = disksfun(z)
u = log(abs(z)) + C;
for j = 1:J
   cj = c(j); u = u + d(j)*log(abs(z-cj));
   for k = 1:N, zck = (z-cj).^(-k); kk = k+(j-1)*N;
      u = u+a(kk)*real(zck)+b(kk)*imag(zck);
   end
   u(abs(z-cj)<=r(j)) = NaN;
end
end

function g = dzdt(t,z)
g = conj(1./z);
for j = 1:J
   zcj = z - c(j); g = g + d(j)./conj(zcj);
   for k = 1:N, kk = k+(j-1)*N;
      g = g - k*(a(kk)+1i*b(kk))./conj(zcj).^(k+1);
   end
end
g = g./abs(g);
end

function [val,isterm,dir] = disksevent(t,z)
dir = zeros(J,1); isterm = ones(J,1); val = abs(z-c.')-r.';
end

end
\end{verbatim}
\smallskip

\begin{verbatim}
function slit1   % Green function exterior to a slit
%
% We seek the function u(z) that is zero on the interval c + r[-1,1] and harmonic outside
% this slit, including at z=infty, except with u(z)~log|z| as z->0.  u is expanded as
%
%   u(z) = log|z| - log|w| + a(1)
%              + SUM_{k=1}^N a(2k)*real(w^(-k)) + a(2k+1)*imag(w^(-k)) 
%
% where w(z) is a Joukowski map of exterior(slit) to exterior(unit disk).
     
c = 3+1i; r = (1-.5i);                        % center and half-length of slit
N = 10; npts = 3*N;                           % no. expansion terms and sample pts
circ = (1+1e-12)*exp(2i*pi*(1:npts)'/npts);
z = c + r*(circ+1./circ)/2;                   % sample points
w = wz(z); rhs = -log(abs(z)) + log(abs(w));  % right-hand side 
A = ones(npts,2*N+1);
for k = 1:N                                   % set up least-squares matrix
   A(:,2*k  ) = real(w.^-k); 
   A(:,2*k+1) = imag(w.^-k); 
end
a = A\rhs;                                    % solve least-squares problem
    
% Contour plot
x = linspace(-5,5,145); y = linspace(-4,4,115); [xx,yy] = meshgrid(x,y);
zz = xx+1i*yy; uu = slit1fun(zz); z = c + r*[-1 1];
plot(z,'b','linewidth',1.4), hold on
levels = -3:.25:-.25; contour(xx,yy,uu,levels,'k'), axis equal, axis([-5 5 -4 4])
set(gca,'xtick',-4:2:4,'ytick',-4:2:4,'fontsize',8)
op = odeset('events',@event);
for t = pi*(1:32)/16, z0 = .01*exp(1i*t); sol = ode23(@dzdt,[0 50],z0,op);
plot(deval(sol,linspace(0,max(sol.x),300)),'k'), end
plot(0,0,'.r','markersize',13), hold off, print -depsc slit1

function u = slit1fun(z)
w = wz(z); u = log(abs(z)) - log(abs(w)) + a(1);
for k = 1:N, u = u + a(2*k)*real(w.^(-k)) + a(2*k+1)*imag(w.^(-k)); end
end

function w = wz(z)
zc = (z-c)/r; sgn = real(zc)>0|(real(zc)==0&imag(zc)>0); sgn = 2*sgn - 1;
w = zc + sgn.*sqrt(zc.^2-1);
end

function g = dzdt(t,z)
zc = (z-c)/r; sgn = real(zc)>0|(real(zc)==0&imag(zc)>0); sgn = 2*sgn - 1;
cw = conj(zc+sgn.*sqrt(zc.^2-1)); dwdzc = conj((1+sgn*zc./sqrt(zc.^2-1))/r);
g = 1./conj(z) - dwdzc./cw;
for k = 1:N, g = g - k*(a(2*k)+1i*a(2*k+1))*dwdzc./cw.^(k+1); end
g = g./abs(g);
end

function [val,isterm,dir] = event(t,z)
dir = 0; isterm = 1; val = abs(wz(z))-1.03;
end

end
\end{verbatim}
\smallskip

\begin{verbatim}
function slits   % Green function exterior to several slits
%
% We seek the function u(z) that is zero on the complex intervals c(j)+r(j)[-1,1]
% j = 1,...,J, and harmonic outside these slits, including at z=infty, except with 
% u(z)~log|z| as z->0.  u is expanded as
%
%   u(z) = log|z| + C + SUM_{j=1}^J {d(j)*log|wj(z)|
%            + SUM_{k=1}^N [a(j,k)*real(wj(z)^(-k)) + b(j,k)*imag(wj(z)^(-k)]},
%
% where wj(z) is a Joukowski map of exterior(slit j) to exterior(unit disk),
% with SUM d(j) = -1; all these coefficients are collected in the vector X.
% The unknowns determined by linear least-squares are C, d(1),...,d(J),
% {a(j,k)}, {b(j,k)}.  A has dimensions 1+J*npts by 1+J*(2N+1).
     
c = [-3 3i 3 2-3i];                          % centers
r = [.5i .5 .3+.3i 1.5+.2i]; J = 4;          % radii
N = 10; npts = 3*N;                          % no. expansion terms and sample pts
circ = (1+1e-8)*exp(2i*pi*(1:npts)'/npts);   % roots of unity
z = []; for j = 1:J
   z = [z; c(j)+r(j)*(circ+1./circ)/2]; end  % sample points on the bndry
A = ones(size(z));                           % constant term
for j = 1:J
   wj = wz(z,j); A = [A log(abs(wj))];       % logarithmic terms
   for k = 1:N                               % set up least-squares matrix
      wck = wj.^(-k);
      A = [A real(wck) imag(wck)];           % algebraic terms
  end
end
A = [A; zeros(1,1+J*(2*N+1))];
A(end,2:(2*N+1):end) = 1;
rhs = [-log(abs(z)); -1];
X = A\rhs;                                   % solve least-squares problem
C = X(1); X(1) = [];                         % extract results
d = X(1:2*N+1:end), X(1:2*N+1:end) = [];
a = X(1:2:end); b = X(2:2:end);

% Contour plot
x = linspace(-5,5,145); y = linspace(-4,4,115); [xx,yy] = meshgrid(x,y);
zz = xx+1i*yy; uu = slitsfun(zz);
for j = 1:J, slit = c(j)+r(j)*[-1 1]; plot(slit+1e-12i,'b','linewidth',1.4), hold on, end
levels = -2.3:.2:-.1; contour(xx,yy,uu,levels,'k'), axis equal, axis([-5 5 -4 4])
set(gca,'xtick',-4:2:4,'ytick',-4:2:4,'fontsize',8), op = odeset('events',@event);
for t = pi*(1:32)/16, z0 = .01*exp(1i*t); sol = ode23(@dzdt,[0 100],z0,op);
plot(deval(sol,linspace(0,max(sol.x),300)),'k'), end
plot(0,0,'.r','markersize',13), hold off, print -depsc slits

function u = slitsfun(z)
u = log(abs(z)) + C;
for j = 1:J
   cj = c(j); w = wz(z,j); u = u + d(j)*log(abs(w));
   for k = 1:N, wk = w.^(-k); kk = k+(j-1)*N;
      u = u+a(kk)*real(wk)+b(kk)*imag(wk);
   end
   u(abs(w)<=1.01) = NaN;
end
end

function w = wz(z,j)
zc = (z-c(j))/r(j); sgn = real(zc)>0|(real(zc)==0&imag(zc)>0); sgn = 2*sgn - 1;
w = zc + sgn.*sqrt(zc.^2-1);
end

function g = dzdt(t,z)
g = 1./conj(z);
for j = 1:J
   zc = (z-c(j))/r(j); sgn = real(zc)>0|(real(zc)==0&imag(zc)>0); sgn = 2*sgn - 1;
   cw = conj(zc+sgn.*sqrt(zc.^2-1)); dwdzc = conj((1+sgn*zc./sqrt(zc.^2-1))/r(j));
   g = g + d(j)*dwdzc./cw;
   for k = 1:N, kk = k+(j-1)*N; g = g - k*(a(kk)+1i*b(kk))*dwdzc./cw.^(k+1); end
end
g = g./abs(g);
end

function [val,isterm,dir] = event(t,z)
dir = zeros(J,1); isterm = ones(J,1);
val = zeros(J,1); for j = 1:J; val(j) = abs(wz(z,j))-1.03; end
end

end
\end{verbatim}
\smallskip

\begin{verbatim}
function cantor  % Green function exterior to a finite approximate Cantor set
%
% We seek the function u(z) that is zero on the real intervals c(j)+r(j)[-1,1]
% j = 1,...,J, and harmonic outside these slits, including at z=infty, except with 
% u(z)~log|z| as z->0.  u is expanded as
%
%   u(z) = log|z| + C + SUM_{j=1}^J {d(j)*log|wj(z)|
%            + SUM_{k=1}^N a(j,k)*real(wj(z)^(-k)),
%
% where wj(z) is a Joukowski map of exterior(slit j) to exterior(unit disk),
% with SUM d(j) = -1; all these coefficients are collected in the vector X.
% The unknowns determined by linear least-squares are C, d(1),...,d(J),
% {a(j,k)}.  A has dimensions 1+J*npts by 1+J*(N+1).
     
m = 5; J=2^m; c = 0; r = 1.5;                 % level of finite Cantor set
for i = 1:m
  r = r/3; c = [c+2*r c-2*r];
end
N = max(2,6-m); npts = round(1.5*N);          % no. expansion terms and sample pts
circ = (1+1e-12)*exp(1i*pi*(.5:npts)'/npts);  % roots of unity
z = []; for j = 1:J
  z = [z; c(j)+r*(circ+1./circ)/2]; end       % sample points on the bndry
A = ones(size(z));                            % constant term
for j = 1:J
  wj = wz(z,j); A = [A log(abs(wj))];         % logarithmic terms
  for k = 1:N, A = [A real(wj.^(-k))]; end    % set up least-squares matrix
end
A = [A; zeros(1,1+J*(N+1))];
A(end,2:N+1:end) = 1;
rhs = [-log(abs(z)); -1];
X = A\rhs;                                    % solve least-squares problem
C = X(1); X(1) = [];                          % extract results
d = X(1:N+1:end); X(1:N+1:end) = []; a = X;

% Contour plot
x = linspace(0,1.8,145); y = linspace(-0.2,0.8,75); [xx,yy] = meshgrid(x,y);
zz = xx+1i*yy; uu = cantorfun(zz);
for j = 1:J, slit = c(j)+r*[-1 1]; plot(slit+1e-12i,'-b','linewidth',1.2), hold on, end
levels = -10.^(-1.2:.1:.2); contour(xx,yy,uu,levels,'k'), axis equal, axis([0 1.8 -.2 .8])
set(gca,'xtick',0:.4:2,'ytick',-.8:.4:.8,'fontsize',5), op = odeset('events',@event);
for t = pi*(1:2:23)/48, z0 = .01*exp(1i*t); sol = ode23(@dzdt,[0 100],z0,op);
z = deval(sol,linspace(0,max(sol.x),300)); plot(z,'k'), plot(conj(z),'k'), end
plot(0,0,'.r','markersize',16), hold off

function u = cantorfun(z)
u = log(abs(z)) + e;
for j = 1:J
  cj = c(j); w = wz(z,j); u = u + d(j)*log(abs(w));
  for k = 1:N, wk = w.^(-k); kk = k+(j-1)*N; u = u+a(kk)*real(wk); end
  u(abs(w)<=1.01) = NaN;
end
end

function w = wz(z,j)
zc = (z-c(j))/r; sgn = real(zc)>0|(real(zc)==0&imag(zc)>0); sgn = 2*sgn - 1;
w = zc + sgn.*sqrt(zc.^2-1);
end

function g = dzdt(t,z)
g = 1./conj(z);
for j = 1:J
  zc = (z-c(j))/r; sgn = real(zc)>0|(real(zc)==0&imag(zc)>0); sgn = 2*sgn - 1;
  cw = conj(zc+sgn.*sqrt(zc.^2-1)); dwdzc = conj((1+sgn*zc./sqrt(zc.^2-1))/r);
  g = g + d(j)*dwdzc./cw;
  for k = 1:N, kk = k+(j-1)*N; g = g - k*a(kk)*dwdzc./cw.^(k+1); end
end
g = g./abs(g);
end

function [val,isterm,dir] = event(t,z)
dir = zeros(J,1); isterm = ones(J,1);
val = zeros(J,1); for j = 1:J; val(j) = abs(wz(z,j))-(1+.01/r); end
end

end
\end{verbatim}
\smallskip

\begin{verbatim}
function bounded   % Harmonic function in a region bounded by disks
%
% We seek the function u with u(z)=0 for |z-c(1)|=r(1)>>1 and u(z)=1 for |z-c(j)|=r(j),
% j = 2,...,J, that is harmonic in the region between these circles.  u is expanded as
%
%   u(z) =  C + SUM_{k=1}^N [a(1,k)*real((z-z(j))^k + b(1,k)*imag((z-z(j))^k]
%             + SUM_{j=2}^J {d(j)*log|z-c(j)|
%             +   SUM_{k=1}^N [a(j,k)*real((z-z(j))^-k + b(j,k)*imag((z-z(j))^-k]}
%
% all these coefficients are collected in the vector X.  The unknowns determined by
% linear least-squares are C, d(2),...,d(J), {a(j,k)}, {b(j,k)}.  A has dimensions
% J*npts by J*(2N+1).
     
c = [0 -1.6+.8i 1.6+.8i];                   % centers
r = [3.8 1 .5]; J = 3;                      % radii
N = 10; npts = 3*N;                         % no. expansion terms and sample pts
circ = exp(2i*pi*(1:npts)'/npts);           % roots of unity
z = []; for j = 1:J
   z = [z; c(j)+r(j)*circ]; end             % sample points on the bndry
A = ones(size(z));                          % constant term
for j = 1:J
   s = 1;
   if j>1
      s = -1; A = [A log(abs(z-c(j)))];     % logarithmic terms
   end  
   for k = 1:N                              % set up least-squares matrix
      zck = (z-c(j)).^(s*k);
      A = [A real(zck) imag(zck)];          % algebraic terms
   end
end
rhs = ones(size(z)); rhs(1:npts) = 0;
X = A\rhs;                                  % solve least-squares problem
C = X(1); X(1) = [];                        % extract results
d = [NaN; X(2*N+1:2*N+1:end)] 
X(2*N+1:2*N+1:end) = [];
a = X(1:2:end); b = X(2:2:end);
    
% Contour plot
x = linspace(-5,5,145); y = linspace(-4,4,115); [xx,yy] = meshgrid(x,y); zz = xx+1i*yy;
uu = boundedfun(zz); z = exp(2i*pi*(0:300)'/150);
plot(c(1)+r(1)*z,'r','linewidth',1.2), hold on
for j = 2:J, disk = c(j)+r(j)*z; fill(real(disk),imag(disk),[.7 .7 1])
plot(disk,'b'), end, levels = .1:.1:.9;
contour(xx,yy,uu,levels,'k'), axis equal, axis([-5 5 -4 4])
set(gca,'xtick',-4:2:4,'ytick',-4:2:4,'fontsize',8), op = odeset('events',@boundedevent);
for t = pi*(1:64)/32, z0 = r(1)*exp(1i*t); sol = ode23(@dzdt,[0 10],z0,op);
plot(deval(sol,linspace(0,max(sol.x),300)),'k'), end
hold off, axis off, print -depsc bounded

function u = boundedfun(z)
u = C*ones(size(z));
for k = 1:N
   zck = (z-c(1)).^k; u = u+a(k)*real(zck)+b(k)*imag(zck);
end
u(abs(z-c(1))>=r(1)) = NaN;
for j = 2:J
   cj = c(j); u = u + d(j)*log(abs(z-cj));
   for k = 1:N
      zck = (z-cj).^(-k); kk = k+(j-1)*N; u = u+a(kk)*real(zck)+b(kk)*imag(zck);
   end
   u(abs(z-cj)<=r(j)) = NaN;
end
end

function g = dzdt(t,z)
g = zeros(size(z));
for k = 1:N, g = g + k*(a(k)+1i*b(k)).*conj(z-c(1)).^(k-1); end
for j = 2:J
   zcj = z - c(j); g = g + d(j)./conj(zcj);
   for k = 1:N, kk = k+(j-1)*N;
      g = g - k*(a(kk)+1i*b(kk))./conj(zcj).^(k+1);
   end
end
g = g./abs(g);
end

function [val,isterm,dir] = boundedevent(t,z)
dir = zeros(J,1); isterm = ones(J,1); val = abs(z-c.')-r.';
end

end
\end{verbatim}
\smallskip

\par}


\begin{thebibliography}{1}

\bibitem{amano}
K. Amano, D. Okano, H. Ogata, and M. Sugihara,
Numerical conformal mappings onto the linear slit domain,
{\em Japan J. Indust.\ Appl.\ Math.} {\bf 29} (2012) 165--186.

\bibitem{axler}
S. Axler, Harmonic functions from a complex analysis viewpoint,
{\em Amer.\ Math.\ Monthly} {\bf 93} (1986) 246--258.


\bibitem{akt}
A. P. Austin, P. Kravanja, and L. N. Trefethen,
Numerical algorithms based on analytic function values at
roots of unity, {\em SIAM J. Numer.\ Anal.} {\bf 52} (2014) 1795--1821.

\bibitem{bb} 
A. H. Barnett and T. Betcke, Stability and convergence
of the method of fundamental solutions for Helmholtz problems on
analytic domains,
{\em J. Comp.\ Phys.} {\bf 227} (2008) 7003--7026.

\bibitem{bergman}
S. Bergman, Functions satisfying certain partial differential
equations of elliptic type and their representation,
{\em Duke Math.\ J.} {\bf 14} (1947) 349--366.

\bibitem{bt}
T. Betcke and L. N. Trefethen, Reviving the
method of particular solutions, {\em SIAM Review} {\bf 47}
(2005) 469--491.

\bibitem{bt2}
T. Betcke and L. N. Trefethen, Computed eigenmodes of
planar regions,
{\em Contemp.\ Math.} {\bf 412} (2006) 297--314.


\bibitem{bourot2}
J. M. Bourot and F. Moreau, Sur l'utilisation de la s\'erie
cellulaire pour le calcul d'\'ecoulements plans de Stokes
en canal ind\'efini: application au cas d'un cylindre circulaire
en translation, {\em Mechanics Research Communications} {\bf 14}
(1987) 187--197.

\bibitem{faraday}
S. J. Chapman, D. P. Hewett, and L. N. Trefethen,
Mathematics of the Faraday cage,
{\em SIAM Rev.} {\bf 57} (2015) 398--417.


\bibitem{cheng}
H. Cheng, On the method of images for systems of
closely spaced conducting spheres,
{\em SIAM J. Appl.\ Math.} {\bf 61} (2001), 1324--1337.

\bibitem{chengg1}
H. Cheng and L. Greengard, On the numerical evaluation of electrostatic
fields in dense random dispersions of cylinders,
{\em J. Comp.\ Phys.} {\bf 136} (1997), 629--639.

\bibitem{chengg2}
H. Cheng and L. Greengard, A method of images for the evlauation of
electrostatic fields in systems of closely spaced conductors,
{\em SIAM J. Appl.\ Math.} {\bf 58} (1998), 122--141.

\bibitem{ckgn}
D. G. Crowdy, E. H. Kropf, C. C. Green, and M. M. S.
Nasser, The Schottky--Klein prime function: a
theoretical and computational tool for applications,
{\em IMA J. Appl.\ Math.} {\bf 81} (2016) 589--628.

\bibitem{crowdy08}
D. G. Crowdy, Geometric function theory: a modern view
of a classical subject, {\em Nonlinearity} {\bf 21}
(2008) T205--T219.

\bibitem{crowdy12}
D. G. Crowdy, Conformal slit maps in applied mathematics,
{\em ANZIAM J.} {\bf 53} (2012) 171--189.

\bibitem{crowdy1}
D. Crowdy and J. Marshall,
Conformal mappings between canonical multiply connected domains,
{\em Comp.\ Meth.\ Funct.\ Th.} {\bf 6} (2006) 59--76.

\bibitem{crowdy}
D. G. Crowdy and J. S. Marshall,
Computing the Schottky--Klein prime function on
the Schottky double of planar domains,
{\em Comp.\ Meth.\ Funct.\ Th.} {\bf 7} (2007) 293--308.

\bibitem{curtiss}
J. H. Curtiss, Interpolation by harmonic polynomials,
{\em J. SIAM} {\bf 10} (1962) 709--736.

\bibitem{ddep}
T. K. De Lillo, T. A. Driscoll, A. R. Elcrat, and
J. A Pfaltzgraff,
Radial and circular slit maps of unbounded
multiply connected circle domains,
{\em Proc.\ Roy.\ Soc.\ A} {\bf 464} (2008) 1719--1737.



\bibitem{dekp}
T. K. De Lillo, A. R. Elcrat, E. H. Kropf and J. A. Pfaltzgraff,
Efficient calculation of Schwarz--Christoffel transformations
for multiply connected domains using Laurent series, 
{\em Comp.\ Meth.\ Funct.\ Th.} {\bf 13} (2013) 307--336.

\bibitem{SC}
T. A. Driscoll and L. N. Trefethen,
{\em Schwarz--Christoffel Mapping,}
Cambridge University Press, 2002.

\bibitem{finn}
M. D. Finn, S. M. Cox, and H. M. Byrne,
Topological chaos in inviscid and viscous
mixers, {\em J. Fluid Mech.} {\bf 493} (2003) 345--361.

\bibitem{fhm}
L. Fox, P. Henrici, and C. Moler,
Approximations and bounds for eigenvalues of elliptic
operators, {\em SIAM J. Numer.\ Anal.} {\bf 4} (1967) 89--102.

\bibitem{gaier}
D. Gaier, {\em Lectures on Complex Approximation,}
Birkh\"auser, 1987.

\bibitem{gander}
W. Gander and J. H\v reb\'i\v cek,
{\em Solving Problems in Scientific Computing
Using Maple and MATLAB,}
4rd ed., Springer, 2004.

\bibitem{gm}
J. B. Garnett and D. E. Marshall,
{\em Harmonic Measure,} Cambridge University Press, 2005.

\bibitem{green}
C. Green, Using the Schottky--Klein prime function
to compute harmonic measure distribution functions of
a class of multiply connected planar domains,
ANZIAM 2018, 8 February 2018.

\bibitem{greensw}
C. Green, M. A. Snipes, and L. A. Ward,
Harmonic measure distribution functions for a
class of multiply connected symmetric slit domains,
manuscript in preparation, March 2018.

\bibitem{greengard}
L. Greengard and V. Rokhlin, A new version of the
fast multipole method for the Laplace equation in three dimensions,
{\em Acta Numer.} {\bf 6} (1997) 229--269.

\bibitem{hack}
W. Hackbusch,
{\em Hierarchical Matrices: Algorithms and Analysis,}
Springer, 2015.

\bibitem{helsing}
J. Helsing and R. Ojala, On the evaluation of layer
potentials close to their sources,
{\em J. Comp.\ Phys.} {\bf 227} (2008) 2899--2921.

\bibitem{henrici}
P. Henrici, A survey of I. N. Vekua's theory of elliptic
partial differential equations with analytic coefficients,
{\em Z. Angew.\ Math.\ Phys.} {\bf 8} (1957) 169--203.

\bibitem{hockney}
R. W. Hockney, A solution of Laplace's equation for a round
hole in a square peg, {\em J. SIAM} {\bf 12} (1964) 1--14.

\bibitem{kk}
L. V. Kantorovich and V. I. Krylov,
{\em Approximate Methods of Higher Analysis,}
Interscience, 1960.

\bibitem{lsn}
J. Liesen, O. S\`ete, and M. Nasser,
Fast and accurate computation of the logarithmic capacity
of compact sets,
{\em Comp.\ Meth.\ Funct.\ Th.} {\bf 17} (2017) 689--713.

\bibitem{malik}
N. H. Malik, A review of the charge simulation
method and its applications,
{\em IEEE Trans.\ Elect.\ Insul.} {\bf 24.1} (1989) 3--20.

\bibitem{moler69}
C. B. Moler, Accurate bounds for the eigenvalues of the
Laplacian and applications to rhombical domains,
Tech.\ Rep.\ CS 121, Dept.\ of Computer Science, Stanford
University, 1969,
{\tt http://i.stanford.edu/TR/CS-TR-69-121.html}.

\bibitem{nasser1}
M. M. S. Nasser, Fast solution of boundary integral equations with
the generalized Neumann kernel,
{\em Elect.\ Trans.\ Numer.\ Anal.} {\bf 44} (2015) 189--229.

\bibitem{nasser}
M. M. S. Nasser, A fast numerical
method for ideal fluid flow in domains with multiple stirrers,
{\em Nonlinearity} {\bf 31} (2018) 815--837.

\bibitem{nassergreen}
M. M. S. Nasser and C. C. Green, 
A fast numerical method for ideal fluid flow in domains
with multiple stirrers,
{\em Nonlinearity} {\bf 31} (2018) 815--387.

\bibitem{nmia}
M. M. S. Nasser, A. H. M. Murid, M. Ismail and E. M. A. Alejaily,
Boundary ingegral equations with the generalized Neumann kernel
for Laplace's equation in multiply connected regions,
{\em Appl.\ Math.\ Comp.} {\bf 217} (2011) 4710--4727.

\bibitem{ogata}
H. Ogata and M. Katsurada,
Convergence of the invariant scheme of the method
of fundamental solutions for two-dimensional
potential problems in a Jordan region,
{\em Japan J. Indust.\ Appl.\ Math.} {\bf 31} (2014) 231--262.



\bibitem{price}
T. J. Price, T. Mullin and J. J. Kobine,
Numerical and experimental characterization of a
family of two-roll-mill flows,
{\em Proc.\ R. Soc.\ Lond.\ A} {\bf 459} (2003) 117--135.

\bibitem{prosnak}
W. J. Prosnak, {\em Computation of Fluid Motions in
Multiply Connected Domains,} G. Braun, Karlsruhe, 1987.

\bibitem{ransford}
T. Ransford, {\em Potential Theory in the Complex Plane,}
Cambridge U. Press, 1995.

\bibitem{ranros}
T. Ransford and J. Rostand, Computation of capacity,
{\em Math.\ Comput.} {\bf 76} (2007) 1499--1520.

\bibitem{rayleigh}
Lord Rayleigh, On the influence of obstacles
arranged in rectangular order upon the properties
of a medium, {\em Philos.\ Mag.} {\bf 32}
(1892) 481--502.

\bibitem{rostand}
J. Rostand, Computing logarithmic capacity with
linear programming, 
{\em Exper.\ Math.} {\bf 6} (1997) 221--238.

\bibitem{rungethm}
C. Runge, Zur Theorie der eindeutigen analytischen Functionen,
{\em Acta Math.} {\bf 6} (1885) 229--244.

\bibitem{runge01}
C. Runge, \"Uber empirische Funktionen und die Interpolation
zwischen \"aquisidtanten Ordinaten,
{\em Z. Math.\ Phys.} {\bf 46} (1901) 224--243.



\bibitem{ssw}
H. Singer, H. Steinbigler, and P. Weiss,
A charge simulation method for the calculation of high voltage
fields, {\em IEEE Trans.\ Pow.\ App.\ Syst.} {\bf 5} (1974)
1660--1668.

\bibitem{tda}
L. N. Trefethen, Ten digit algorithms, Numer.\ Anal.\
Rep.\ 05/13, Oxford U. Computing Lab., 2005,
\verb|https://people.maths.ox.ac.uk/trefethen/publication/PDF/2005_114.pdf|.

\bibitem{atap}
L. N. Trefethen, {\em Approximation Theory and
Approximation Practice,} SIAM, 2013.



\bibitem{walsh29}
J. L. Walsh, The approximation of harmonic functions
by harmonic polynomials and by harmonic rational
functions, {\em Bull.\ Amer.\ Math.\ Soc.} {\bf 35} (1929)
499--544.


\bibitem{zaviska}
F Z\'avi\v ska, \"Uber die Beugung elektromagnetischer Wellen
an parallelen, unendlich langen Kreiszylindern,
{\em Ann.\ Phys.} {\bf 345} (1913) 1023--1056.

\end{thebibliography}
\end{document}